\documentclass{amsart}
\usepackage{enumerate,bm}
\usepackage{graphicx}
\usepackage[outdir=./]{epstopdf}
\usepackage{hyperref}
\usepackage{xcolor}
\usepackage{soul}
\pagecolor{white}
\usepackage[linesnumbered,ruled]{algorithm2e}
\usepackage{multirow}
\usepackage{amssymb,amsmath}
\usepackage{tensor}
\usepackage{cancel}
\usepackage{mdframed}
\usepackage{algorithmic}
\usepackage{makecell}
\usepackage{booktabs}
\usepackage{subcaption}

\let\oldnl\nl
\newcommand{\nonl}{\renewcommand{\nl}{\let\nl\oldnl}}

\newtheorem{prop}{Proposition}

\newtheorem{coro}{Corollary}
\newcommand{\R}{\mathbb R}
\newcommand\ChangeRT[1]{\noalign{\hrule height #1}}

\title[Low Precision Projected Iterated Tikhonov]{Projected iterated Tikhonov regularization in low precision}

\author{C. Drum}
\author{J.~G. Nagy}
\author{L. Onisk}

\address{Department of Mathematics, Emory University, Atlanta, GA 30322}
\email{cndrum@emory.edu}
\email{jnagy@emory.edu}
\email{lonisk@emory.edu}

\begin{document}

\begin{abstract}
We investigate the regularizing behavior of an iterative Krylov subspace method for the solution of linear inverse problems in precisions lower than double precision. Recent works have considered the projection of iterated Tikhonov methods using Krylov subspaces for both computational efficiency and an additional regularizing effect. To investigate the regularizing behavior of this projected algorithm applied to problems that are naturally severely ill-posed, we formulate the iterates as a filtered solution using the preconditioned Landweber method with a Tikhonov-type preconditioner in a Krylov subspace. Through numerical examples simulating multiple low precision choices, we showcase the filtering properties of the method and the achievement of comparable working accuracy applied to discrete inverse problems (i.e., to within a few decimal places in relative error) compared to results computed in traditional double precision.

\end{abstract}

\maketitle

\section{Introduction} \label{sec:intro}
Herein, we are interested in computing approximate solutions to linear least-squares problems of the form
\begin{equation} \label{minProb}
    \min_{x\in \mathbb{R}^{n}} \left\|Ax-b\right\|,
\end{equation}
in computing environments that support low precision and where data is often stored with less significant digits, e.g., medical tomography \cite{Siltanen22}, industrial applications \cite{Boas2012}, and training for deep learning \cite{Ruthotto25}. Here, and throughout we denote the Euclidean norm by $\|\cdot\|$. For some of these applications the singular values of the matrix $A \in \mathbb{R}^{m \times n}$ can decay without significant gap and cluster at the origin (i.e., the matrix is ill-conditioned). These types of problems are commonly referred to as linear discrete ill-posed problems. Classically, they can arise through the discretization of Fredholm integral equations of the first kind \cite{engl1996regularization, Ha98}, but can also arise in modern large scale data problems such as, for example, the training of the random feature model in machine learning \cite{RB07}.

In applications, the vector $b \in \mathbb{R}^m$ in \eqref{minProb} often represents observed measurements that are corrupted by error arising through noise contamination or instrumentation disturbances. Numerically, this error may also include truncation error or approximation error. We denote this error by $e$ so that
\begin{equation*}
    b = b_{true} + e
\end{equation*}
where $b_{true}$ is the unknown error-free vector associated with $b$. The goal when considering linear discrete ill-posed problems is to determine the solution of minimal Euclidean norm, $x^{\dagger}$, of the unavailable least-squares problem 
\begin{equation} \label{trueMinProb}
   x^{\dagger}= \arg\min_{x\in \mathbb{R}^{n}} \left\|Ax-\hat{b}\right\|.
\end{equation}
In the linear discrete inverse problem literature it's well understood that because $b$ is contaminated by error and the singular values of $A$ cluster at the origin that the solution of \eqref{minProb} of minimal Euclidean norm is usually a poor approximation to $x^{\dagger}$. The common strategy is to regularize the problem - that is, to replace the problem \eqref{minProb} with a nearby problem that is less sensitive to the error in $b$.

Recent progress in GPU hardware that supports floating point computations in mixed-to-low precision has enabled the numerical linear algebra community to develop robust decomposition strategies as well as accelerated iterative methods; see the surveys \cite{MixedPrecisionSurvey2021, Higham22}. Iterative methods that utilize Krylov subspaces using one or more precisions lower than fp64 (i.e., double precision) have recently been considered in the literature for well-posed problems, see e.g., \cite{carson2025mixed,CarsonHighamPranesh20, Carson2}. While some works including \cite{Renaut25} and \cite{deSturler25} have considered low or mixed precision in the inverse problem community, comparatively less has been studied, especially for iterative Krylov methods. 

Recent work by Li in \cite{Haibo24} considered LSQR in low precision for linear inverse problems with the conclusion that `LSQR does not need to be solved in high precision'. However, the standard practice in many iterative Krylov methods today is to utilize a hybrid approach whereby one both projects the problem into a subspace and then performs regularization through a variational approach, rather than solely regularizing implicitly via subspace projection; see the review \cite{ChungGazzola_survey}. In a similar vein to the `project then regularize' strategy just mentioned, methods based on projected iterated Tikhonov regularization with a nonstationary regularization parameter using a Krylov subspace have garnered recent attention; see \cite{bianchi2025convergence,BGOPR25,BOR23,chang2024relaxed}.

Classically, iterated Tikhonov regularization is used to approximate the solution to \eqref{minProb} and at its $k^{th}$ iterate is given by 
\begin{equation} \label{iterTik.}
    x^{(k)} = x^{(k-1)} + \left(A^TA+\alpha^2I\right)^{-1}A^Tr^{(k-1)}
\end{equation}
with the residual denoted by $r^{(k-1)} = b-A^{(k-1)}$; see \cite{buccini2017iterated, HG98}. Here, $\alpha>0$ is the regularization parameter that balances the fidelity of the approximation \eqref{iterTik.} to \eqref{minProb}. A major drawback to \eqref{iterTik.} is that if $A$ is large then a large system of equations with the matrix $A^TA+\alpha^2I$ has to be solved at each iteration. In \cite{BOR23} a low rank approximation to $A$ at each step $k$ using the Arnoldi process was used as an efficient alternative to the work by Donatelli and Hanke in \cite{DH13} which replaced $A$ with an approximation $C$ in \eqref{iterTik.} whose structure allows for fast computation (e.g., using the discrete Fourier transform where appropriate). To generalize the work of \cite{BOR23}, the authors of \cite{BGOPR25} additionally considered the use of Golub-Kahan bidiagonalization (GKB) when both $A$ and $A^T$ are available and $A$ is possibly non-square. Additionally, they investigated use of a precomputed GKB subspace, which was found to provide better accuracy in practice then the technique of iteratively increasing the subspace in \cite{BOR23}. This method is referred to as the \emph{projected iterated Tikhonov} (PIT) method.

Herein, we investigate the regularizing behavior of the PIT method applied to \eqref{minProb} in low precision with a nonstationary regularization parameter scheme. To do so, we formulate the iterates of PIT as a projected preconditioned Landweber method with a Tikhonov-type preconditioner that may be computed in a precision lower than fp64 using the GKB process. When considering a mixed precision approach to the Tikhonov problem in \cite{NO24}, it was shown that the precision of the Tikhonov-type preconditioner used in the Landweber method controlled the efficacy of the solution. Alternatively, here we demonstrate that the efficacy of the solution in $\mathbb{R}^n$ can be strongly dependent on the quality of the Krylov basis vectors which, in turn can depend on the working precision. Numerically, we illustrate that when computed in sufficiently low precision that the GKB process requires reorthogonalization to retain the same solution quality as compared to when the method is computed entirely in fp64.

While works including \cite{buccini2017iterated,BGOPR25,BOR23,DH13} consider a nonstationary regularization parameter scheme utilizing a modified Levenberg-Marquardt iteration to solve a damped version of the discrepancy principle, we consider an alternative method based on a secant-update method using the discrepancy principle which requires fewer user chosen parameters and was first proposed in \cite{gazzola_novati} for an iterative Krylov subspace method for linear discrete inverse problems. To gain insight into the effect that the nonstationary regularization parameters have on the PIT method, we derive both the stationary and nonstationary projected filter factors and confirm numerically that the nonstationary filter factors match those that are experimentally realized. We additionally demonstrate in our numerical results that PIT in low precision can give comparable working accuracy (i.e., to within a few decimal places in relative error) against results computed in fp64, which supports its use in modern applications that natively support low precision floating-point arithmetic.

The rest of the paper is structured as follows. In Section \ref{sec:background} we provide necessary background on topics including the Landweber iteration, its connection to iterated Tikhonov regularization, and filtering methods. Section \ref{sec2} introduces necessary background for constructing the PIT method with a fixed regularization parameter and its projected filtering properties before moving onto the nonstationary PIT method in Section \ref{sec3} and its filtering properties. Numerical results and relevant preliminaries are provided in Section \ref{sec:results} and concluding remarks are given in Section \ref{sec:conc}.

\section{Background} \label{sec:background}
This section provides necessary background and notation for representing the regularized solution of a linear discrete ill-posed problem as a filtered solution. Additionally, background on the Landweber method and its preconditoned variant are provided, as well as its connection to iterated Tikhonov regularization.

\subsection{Filtered Solution of Linear Discrete Ill-posed Problems}

The approximate solution of \eqref{minProb} for inverse problems has been well studied using the singular value decomposition (SVD) of $A=U\Sigma V^T$, where $U \in \mathbb{R}^{m \times m}$ and $V \in \mathbb{R}^{n \times n}$ are orthogonal matrices whose columns represent the left and right singular vectors of $A$, respectively. The diagonal matrix $\Sigma \in \mathbb{R}^{m \times n}$ contains the singular values given by $\sigma_j$ for $j=1,2,\dots,n$ on its main diagonal ordered in non-increasing fashion. Using the SVD, one can write the least-squares solution of \eqref{minProb} as a linear combination of the right singular vectors with the $j^{th}$ coefficient given by $u_j^Tb/\sigma_j$. Because the singular values of $A$ cluster near numerical zero, the error $e$ contaminating $b$ is propagated into the computed solution for high index values $j$ corresponding to high frequency information \cite{Ha98, HNO}.

To prevent this propagation of error when approximating the solution of \eqref{minProb}, the conditioning of the problem can be improved by employing filtering methods which can be formulated as a modified inverse solution \cite{HNO}. A filtered solution is of the form
\begin{equation} \label{filteredSoln}
    x_{filt.} = \sum_{j=1}^n  \phi_j \frac{u_j^Tb}{\sigma_j}v_j
\end{equation}
where an intelligent selection of the filter factors $\phi_j \in \mathbb{R}$ can remove deleterious components of the approximate solution corresponding to high frequency information that is dominated by noise. This may be considered a form of regularization. To achieve a meaningful solution, $\phi_j$ are chosen to be $\approx 1$ for small $j$ and  $\approx 0$ for large $j$.

\subsection{Iterated Tikhonov Regularization as Preconditioned Landweber}

Here, we provide an overview of the Landweber iteration and its preconditioned variant - iterated iterated Tikhonov regularization for the solution of linear discrete ill-posed problems. 

\subsubsection{The Landweber Iteration} \label{sec:Landweber}
An alternative strategy to obtain a regularized solution of \eqref{minProb} is to apply an iterative method and terminate the iterations before error is propagated. The discrete Landweber iteration \cite{Landweber} is one such method whose $k^{th}$ iterate is given by
\begin{equation} \label{landweber}
    x^{(k)} = x^{(k-1)} + \zeta A^T\left(b-Ax^{(k-1)}\right)
\end{equation}
where $\zeta \in \left(0,1/\sigma_1^2\right]$ is a relaxation term with $\sigma_1$ the largest singular value of $A$. Without loss of generality herein, we choose $\zeta=1$ since $A$ can be scaled appropriately; see \cite{engl1996regularization} for further discussion. In practice, Landweber is well known to have slow convergence towards a meaningful solution. To understand this we may consider the $k^{th}$ iterate viewed through the lens of its filter factors \cite{HankeHansen93}:
\begin{equation} \label{filterLand}
        x^{(k)} = \sum_{j=1}^{n} \Big(1-(1-\sigma_j^2)^{k} \Big) \frac{u_j^Tb}{\sigma_j}v_j
\end{equation}
where $\phi_j^{(k)} = \left(1-(1-\sigma_j^2)^{k} \right) \in \mathbb{R}$ denotes the $j^{th}$ filter factor a the $k^{th}$ iteration. 

When $A$ is a scaled matrix whose largest singular value is one, only the first filter factor of the $k^{th}$ iteration, $\phi^{(k)}_1$, will be $1$, with subsequent filter factors of the same iteration rapidly decaying to zero. When $k$ increases, the early filter factors, i.e., $\phi^{(k)}_j$ for small $j$, cluster near $1$. Because of this behavior, it can take many iterations to fully capture the dominant right singular vectors that constitute a meaningful approximate solution for inverse problems. Employing a preconditioner is often utilized to alleviate this issue.

\subsubsection{Preconditioning Landweber with a Tikhonov-type preconditioner} \label{sec:preconLand}
Here, we consider preconditioning the Landweber iteration \eqref{landweber} with a Tikhonov-type preconditioner. We clarify that there are various possible preconditioning strategies other than Tikhonov that could be used. It should be understood by our usage of the wording `preconditioned Landweber' in this work that we only consider Tikhonov-type preconditioners; see \cite{Estatico08,Bertero97} for further discussions and other preconditioner types that could be considered for the solution of linear discrete ill-posed problems.

We consider the right-preconditioned least-squares problem
\begin{equation} \label{preconLS}
    \min_{x\in \mathbb{R}^{n}} \left\|AM^{-1}Mx-b\right\|
\end{equation}
with nonsingular preconditioner $M \in \mathbb{R}^{n \times n}$. Denote $\hat{A}=AM^{-1}$ and $\hat{x}=Mx$. Then, we may use Landweber \eqref{landweber} with $\zeta=1$ to solve \eqref{preconLS}, which may be simplified to
\begin{equation} \label{preconLand}
        x^{(k)} = x^{(k-1)} + \left(M^TM\right)^{-1}A^T\left(b-Ax^{(k-1)}\right).
\end{equation}
We note that when $M^TM = \left(A^TA + \alpha^2 I\right)$, that the preconditioned Landweber method is equivalent to iterated Tikhonov regularization \eqref{iterTik.}. Using the SVD of $A$, we may write
\begin{equation*}
M^TM = \left(A^TA + \alpha^2 I\right) = \left(V^T\Sigma^T\Sigma V + \alpha^2I\right) = V^TD^2V
\end{equation*}
where $D^2 = \left(\Sigma^2 + \alpha^2I\right)$ and $\Sigma^2 = \Sigma^T\Sigma \in \mathbb{R}^{n \times n}$. Thus, one may consider the preconditioning matrix $M$ of \eqref{preconLS} as $M = \left(\Sigma^2 + \alpha^2I\right)^{1/2}V^T=DV^T$: the product of a diagonal matrix with the right singular matrix of $A$. We will use this general formulation when considering the filtered solution of the PIT method.

\section{Projected Iterated Tikhonov via Golub-Kahan Bidiagonalization} \label{sec2}

We begin this section with a brief review of the Golub-Kahan bidiagonalization process for forming a Krylov subspace. We then construct the PIT method with a fixed regularization parameter by preconditioning a projected Landweber method. Where possible we reduce possible exposition where it may be inferred from the same general process in Section \ref{sec:preconLand}. To better understand the regularizing behavior of PIT with a fixed regularization parameter, we derive its filter factors by writing it as a filtering method.

\subsection{The Golub-Kahan Bidiagonalization process}

When matrix-vector products with both $A$ and $A^T$ are available, application of $1\leq p\ll n$ iterations of GKB to $A$ with initial vector $b\ne 0$ generically gives the pair of decompositions
\begin{equation}\label{gkb}
AV_p = U_pB_p \quad \& \quad A^T U_{p-1} = V_p B_{p,p}^T,
\end{equation}
where the matrices $U_p\in{\R}^{m\times(p+1)}$ and $V_p\in{\R}^{n\times p}$ have orthonormal columns, $U_{p-1}\in{\R}^{m\times p}$ consists of the first $p$ columns of $U_{p}$, and $U_{p}e_1=b/\|b\|$.
Here, and throughout this paper $e_1=[1,0,\ldots,0]^T$ denotes the first principal axis vector whose size should be apparent from context. The range of $V_p$, denoted by ${\mathcal R}(V_p)$, is the Krylov subspace
\begin{equation*}\label{krylov}
{\mathbb K}_p\left(A^TA,A^Tb\right)=
{\rm span}\left\{A^Tb,\left(A^TA\right)A^Tb,\ldots,(A^TA)^{p-1}A^Tb\right\}.
\end{equation*}
The matrix $B_p \in \mathbb{R}^{(p+1) \times p}$ is lower bidiagonal with positive subdiagonal entries and the matrix $B_{p,p}\in{\R}^{p\times p}$ denotes its leading $p \times p$ submatrix. We assume here that $p$ is small enough so that the decompositions \eqref{gkb} with the stated properties exist. This is the generic situation. The dominating computational effort required to determine the decompositions \eqref{gkb} is the sequential evaluation of $p$ matrix-vector products with the matrices $A$ and $A^T$; see, e.g., \cite[Section 10.4.4]{GVL} for an algorithm. 

We mention that for sufficiently large problems and when computations in low enough floating point precision are done that issues relating to both loss of orthogonality amongst the Krylov basis vectors and overflow or underflow of computations can occur. To the former, in Section \ref{sec:experiments} we describe how reorthogonalization can be necessary to retain orthonormality of the basis vectors $V_p$ of GKB. For a discussion on GKB and orthogonality considerations see \cite{paige1982lsqr}. While we did not experience any of the issues to the ladder, we point the interested reader to \cite{brown2025inner, brown2025h} which consider the solution of \eqref{minProb} using iterative subspace methods that avoid using inner products.

\subsection{PIT with fixed $\alpha$} \label{sec:stationaryPIT}

After $p$ steps of GKB we may write the projected least-squares problem of \eqref{minProb} using \eqref{gkb} as
\begin{equation} \label{projected_minProb}
    \min_{y\in \mathbb{R}^{p}}\left\|B_py- \tilde{b} \right\|
\end{equation}
where $\tilde{b} = U^T_pb \in \mathbb{R}^{p+1}$. If we consider the right-preconditioned projected least-squares problem in the same vein as \eqref{preconLS} by
\begin{equation} \label{preconditioned_projected_minProb}
    \min_{y\in \mathbb{R}^{p}}\left\|B_p\tilde{M}^{-1}\tilde{M}y- \tilde{b} \right\|,
\end{equation}
then we may again use Landweber \eqref{landweber} to solve \eqref{preconditioned_projected_minProb} which may be simplified to 
\begin{equation} \label{preconProj_Land}
    y^{(k)} = y^{(k-1)} + \left(B_p^TB_p + \alpha^2 I\right)^{-1}B^T_p\left(\tilde{b}-B_py^{(k-1)}\right).
\end{equation}
We note that an approximate solution to \eqref{projected_minProb} can be projected back to $\mathbb{R}^n$ by application of $V_p$ from GKB, i.e., $x^{(k)} = V_py^{(k)}$. As such, the $k^{th}$ iterate of PIT may be represented by projecting each update step from projected preconditioned Landweber back to $\mathbb{R}^n$ as follows
\begin{equation} \label{PIT}
\begin{split}
    x^{(k)} &= x^{(k-1)} + V_p\left(B_p^TB_p + \alpha^2 I\right)^{-1}B^T_p\left(\tilde{b}-B_py^{(k-1)}\right)\\
    &= x^{(k-1)} + V_p\left(B_p^TB_p + \alpha^2 I\right)^{-1}B^T_p\tilde{r}^{(k-1)}.
    \end{split}
\end{equation}
In practice, the $k^{th}$ projected iterate, $y^{(k)}$, is only projected back to $\mathbb{R}^n$ when the algorithm terminates since the projected residual retains the same norm value, i.e., $\|b-Ax^{(k)}\| = \|\tilde{b}-B_py^{(k)}\|$, which is used for determination of the nonstationary regularization parameter as well as the stopping criterion (see Sections \ref{sec:regPar} and \ref{sec:nonStationaryPIT}). In what follows we minimize explicit reference to projection into the Krylov subspace using the GKB process, as this should be inferred throughout the remainder of the work when referring to PIT. 

\subsection{Filtered Solution of PIT with fixed $\alpha$} \label{sec:filteredPIT}
We now turn to the consideration of the filtering behavior of PIT starting with a fixed regularization parameter. Following from the previous section that $\tilde{M} \in \mathbb{R}^{p \times p}$ is a nonsingular preconditioner for the projected least-squares problem \eqref{projected_minProb}, we represent it by first defining the SVD $B_p = \hat{U}\hat{\Sigma}\hat{V}^T$ where $\hat{U} \in \mathbb{R}^{(p+1)\times (p+1)}$ and $\hat{V} \in \mathbb{R}^{p\times p}$ are orthogonal matrices, and $\hat{\Sigma} \in \mathbb{R}^{(p+1) \times p}$ is diagonal so that
\begin{equation}
        \tilde{M}^T\tilde{M} = \left(B_p^TB_p + \alpha^2I\right) =  \hat{V}\left(\hat{\Sigma}^2 + \alpha^2I\right)\hat{V}^T = \hat{V}\hat{D}^2\hat{V}^T
\end{equation}
where $\hat{\Sigma}^2 = \hat{\Sigma}^T\hat{\Sigma} \in \mathbb{R}^{p \times p}$ diagonal. As such, we may represent our preconditioner as $\tilde{M} = \left(\hat{\Sigma}^2 + \alpha^2I\right)^{1/2}\hat{V}^T = \hat{D}\hat{V}^T$. We now present the PIT method with fixed $\alpha$ as a filtering method.

\begin{prop}{(Filtered Solution of Stationary PIT)\\} \label{prop1}
After $p$ steps of GKB, the $k^{th}$ projected iterate, $y^{(k)}$, of PIT given by \eqref{preconProj_Land} with preconditioner $\tilde{M}$ and fixed regularization parameter $\alpha >0$ may be written as 
\begin{equation} \label{PITA_projSoln}
	y^{(k)} = \sum_{i=1}^p \psi^{(k)}_i\frac{\hat{u}_{i}^T\tilde{b}}{\hat{\sigma}_i}\hat{v}_{i}
\end{equation}
with filter factors per $k^{th}$ iterative step given by $\psi^{(k)}_i= 1 - \left(\frac{\alpha^2}{\hat{\sigma}_i^2+\alpha^2}\right)^{k}$ for $i=1,2,\dots,p$.
\begin{proof}
\color{black}
We begin by writing the preconditioned least-squares problem \eqref{preconditioned_projected_minProb} as
$$ \min_{y\in \mathbb{R}^{p}}\left\|B_p\tilde{M}^{-1}\tilde{M}y- \tilde{b} \right\| =  \min_{y\in \mathbb{R}^{p}}\left\|\hat{B}_p\hat{y}- \tilde{b} \right\|.$$
whose $k^{th}$ solution by application of Landweber is given by
\begin{equation*}
	\hat{y}^{(k)} = \hat{y}^{(k-1)} + \hat{B}_{p}^T\left(\tilde{b}-\hat{B}_{p}\hat{y}^{(k-1)}\right) 
\end{equation*}
where we have used that $\hat{B}_p = B_{p}\tilde{M}^{-1}$ and $\hat{y}=\tilde{M}y$. Using our definition $\tilde{M} = \left(\hat{\Sigma}^2 + \alpha^2I\right)^{1/2}\hat{V}^T = \hat{D}\hat{V}^T \in \mathbb{R}^{p \times p}$ and the SVD of $B_{p}$, we may rewrite $\hat{B}_{p}$ as follows
\begin{equation*}
		\hat{B}_{p} = B_{p}\tilde{M}^{-1} = B_{p}\hat{V}\hat{D}^{-1}
		= \hat{U}\hat{\Sigma}\hat{V}^T\hat{V}\hat{D}^{-1}.
\end{equation*}
We may then write the following
\begin{equation*}
\hat{B}_{p}^T\hat{B}_{p} = \hat{D}^{-1} \hat{V}^T\hat{V}\hat{\Sigma}\hat{U}^T\hat{U}\hat{\Sigma}\hat{V}^T\hat{V}\hat{D}^{-1} \in \mathbb{R}^{p \times p}
\end{equation*}
so that 
\begin{equation} \label{prop2deriv}
	\begin{split}
    		\hat{y}^{(k)} &= \sum_{j=0}^{k-1} \left(I - \hat{B}_{p}^T\hat{B}_{p}\right)^j\hat{B}_{p}^T\tilde{b} \\
		&= \sum_{j=0}^{k-1} \left(I - \hat{D}^{-1}\hat{V}^T\hat{V}\hat{\Sigma}^T\hat{\Sigma}\hat{V}^T\hat{V}\hat{D}^{-1}\right)^j\hat{D}^{-1}\hat{V}^T\hat{V}\hat{\Sigma}^T\hat{U}^T\tilde{b}.
	\end{split}
\end{equation}
Under orthogonality, \eqref{prop2deriv} simplifies to
\begin{equation*}
	\hat{y}^{(k)} =  \sum_{j=0}^{k-1} \left(I - \hat{D}^{-1}\hat{\Sigma}^2\hat{D}^{-1}\right)^j\hat{D}^{-1}\hat{\Sigma}^T\hat{U}^T\tilde{b}.
\end{equation*}
Denoting $\hat{D}^{k-1}=\sum_{j=0}^{k-1} \left(I - \hat{D}^{-1}\hat{\Sigma}^2\hat{D}^{-1}\right)^j$, we may write in terms of the $k^{th}$ iterative solution, $y^{(k)}$, by using that $y^{(k)} = \tilde{M}^{-1}\hat{y}^{(k)}$ so that
\begin{equation*}
	\begin{split}
		y^{(k)} &= \tilde{M}^{-1}\hat{y}^{(k)}
		= \hat{V}\hat{D}^{-1}\hat{D}^{k-1}\hat{D}^{-1}\hat{\Sigma}^T\hat{U}^T\tilde{b} \\
		&= \sum_{i=1}^p \hat{d}_i^{k-1}\frac{\hat{\sigma}_i}{\hat{d}_{i}^2}\left(\hat{u}_{i}^T\tilde{b}\right)\hat{v}_{i}
	\end{split}
\end{equation*}
where $\hat{d}_i^{k-1} \in \mathbb{R}$ denotes the $i^{th}$ diagonal entry of $\hat{D}^{k-1}$ and may be written as a geometric series
\begin{equation*}
		\hat{d}_i^{k-1} = \sum_{j=0}^{k-1} \left(1 - \frac{\hat{\sigma}_i^2}{\hat{d}_{i}^2}\right)^j
		= \frac{1 - \left(1 - \frac{\hat{\sigma}_i^2}{\hat{d}_{i}^2}\right)^{k}}{\hat{\sigma}_i^2/{\hat{d}_{i}^2}}.
\end{equation*}
With this we may write $y^{(k)}$ as
\begin{equation*}
	\begin{split}
		y^{(k)} &= \sum_{i=1}^p \hat{d}_i^{k-1}\frac{\hat{\sigma}_i}{\hat{d}_{i}^2}\left(\hat{u}_{i}^T\tilde{b}\right)\hat{v}_{i} \\
		&= \sum_{i=1}^p \left(1 - \left(1 - \frac{\hat{\sigma}_i^2}{\hat{d}_{i}^2}\right)^{k}\right)\frac{\hat{u}_{i}^T\tilde{b}}{\hat{\sigma}_i}\hat{v}_{i} = \sum_{i=1}^p \left(1 - \left(\frac{\alpha^2}{\hat{\sigma}_i^2 + \alpha^2}\right)^k\right)\frac{\hat{u}_{i}^T\tilde{b}}{\hat{\sigma}_i}\hat{v}_{i}.
	\end{split}
\end{equation*}
\end{proof}
\end{prop}

With this result, we point out that Proposition \ref{prop1} is the projected analog of the filtered solution of iterated Tikhonov regularization; see \cite{Ha98, NO24}. We comment here that the filtering behavior of \eqref{PITA_projSoln} is improved versus the general behavior of \eqref{filterLand} since for the same matrix $B_p$ and corresponding singular values $\hat{\sigma}_i$ we will have that
$$\left(\frac{\alpha^2}{\hat{\sigma}_i^2+\alpha^2}\right)^k<(1-\hat{\sigma}_i^2)^{k} $$
for $\alpha>0$. Depending on $\alpha$ and $k$, the filter factors of \eqref{PITA_projSoln} will be closer to $1$ for smaller indices $i$ compared with the non-preconditioned variant \eqref{filterLand} which is preferred as mentioned in Section \ref{sec:Landweber} for Landweber style iterative methods. This filtering process may be tuned not only through the iterate $k$, but also by iteratively updating the regularization parameter at each step using the projected residual and knowledge of the noise contaminating the right-hand side of the linear system resulting in what is often referred to as an \emph{nonstationary} iterative regularization method.

\section{Nonstationary PIT in Low Precision} \label{sec3}

We begin this section by introducing the secant-update method which determines our regularization parameters for the PIT method based on the discrepancy principle. The PIT method with nonstationary regularization parameters is then presented as well as our stopping criterion. We conclude by considering the filtering of the nonstationary PIT method.

\subsection{Nonstationary Regularization Parameter Determination} \label{sec:regPar}

Here, we propose the strategy for determining the regularization parameters for the PIT method which we will refer to as the secant-update method or just the secant method originally introduced by Gazzola and Novati in \cite{gazzola_novati}. Interpreted geometrically as a zero finder, the approach is to enforce a nonlinear discrepancy condition at each iteration and update the regularization parameter by a step of a secant method. The formulation of the method in \cite{gazzola_novati} is based on the Arnoldi process; see e.g., \cite{Saad} for more details. Here, we reformulate the scheme using the GKB process so as to integrate with the PIT method with fixed $p$. In doing so, we discuss the discrepancy principle as a means to determine the regularization parameter at each step of PIT.

We assume to know, or have a fairly accurate estimate of an upper bound for the norm of the noise contaminating the available linear system given by $\| b - b_{true} \| = \|e\| \leq \varepsilon$. At the $k^{th}$ iterate of an iterative method, the \emph{discrepancy principle} prescribes that $\alpha>0$ be chosen so that
\begin{equation} \label{DP}
\left\|Ax^{(k)} - b\right\|=\eta\varepsilon
\end{equation}
where $\eta>1$ is a user-specified constant that is independent of $\varepsilon$ and whose value may be interpreted as how much one `trusts' the upper bound on the error contaminating $b$. Here, using \eqref{PIT} we may interpret that $x^{(k)}$ is a function of $\alpha$. We mention that there are other strategies to determine $\alpha$ including the L-curve criterion and generalized cross validation; see e.g., \cite{ChungGazzola_survey, HNO, RR13} for further discussion.

As discussed in Section \ref{sec:stationaryPIT}, the $k^{th}$ projected iterate of PIT takes the form 
\begin{equation*}
y^{(k)}(\alpha) = y^{(k-1)}(\alpha) + \left(B_p^TB_p + \alpha^2 I\right)^{-1}B^T_p\left(\tilde{b}-B_py^{(k-1)}(\alpha)\right)
\end{equation*}
where the notation $y^{(k)}(\alpha)$ is used to denote that the $k^{th}$ projected solution is a function of $\alpha$. Recalling that $x^{(k)} = V_p y^{(k)}$, the discrepancy at each iteration is defined by
\begin{equation*}\label{eq:secantdisc}
\phi^{(k)}(\alpha)
= \left\| b - A x^{(k)}\left(\alpha\right) \right\|
= \left\| \tilde{b} - B_{p} y^{(k)}(\alpha) \right\|,
\end{equation*}
and therefore the discrepancy principle is satisfied at the $k^{th}$ iterate when
\begin{equation}\label{eq:secantDP}
\phi^{(k)}(\alpha) \leq \eta \varepsilon.
\end{equation}
The secant-update approach consists of finding $\alpha>0$ such that nonlinear equation \eqref{eq:secantDP} with equality is satisfied. Following the same approach as in \cite{gazzola_novati}, we employ a linear approximation of the discrepancy:
\begin{equation}\label{eq:secantlinear}
\phi^{(k)}(\alpha) \approx \gamma^{(k)} + \alpha \beta^{(k)},
\end{equation}
where 
\begin{equation*}
\gamma^{(k)} = \phi^{(k)}(0)
= \left\| \tilde{b} - B_{p} y^{(k)}(0) \right\|
\end{equation*}
is the norm of the residual corresponding to the unregularized solution, and where
\begin{equation*}
\beta^{(k)} \approx 
\frac{\phi^{(k)}\left(\alpha^{(k-1)}\right) - \gamma^{(k)}}{\alpha^{(k-1)}},
\end{equation*}
is the discrepancy with corresponding $\alpha$ computed during the previous iteration. Imposing \eqref{eq:secantDP} with equality, i.e.,
\begin{equation*}
\phi^{(k)}(\alpha^{(k)}) = \eta \varepsilon,
\end{equation*}
and using \eqref{eq:secantlinear}, we obtain the update
\begin{equation*}
\alpha^{(k)}
= \frac{\eta\varepsilon - \gamma^{(k)}}
       {\phi^{(k)}\left(\alpha^{(k-1)}\right) - \gamma^{(k)}}\, \alpha^{(k-1)}.
\end{equation*}
This expression can be interpreted as a step of a secant method applied to $\phi^{(k)}(\alpha)$, with the point $(0,\gamma^{(k)})$ as the left 
interpolation node. The authors of \cite{gazzola_novati} noted that when the size of the Krylov subspace is small (i.e., $p$ is small) that $\gamma^{(k)} > \eta\varepsilon$ can occur, which can cause instability. In such situations we may instead use
\begin{equation*}\label{eq:secantupdate}
\alpha^{(k)} 
= \left| 
     \frac{\eta\varepsilon - \gamma^{(k)}}
          {\phi^{(k)}\left(\alpha^{(k-1)}\right) - \gamma^{(k)}}
   \right|\alpha^{(k-1)}.
\end{equation*}
The secant-update method is provided in Algorithm \ref{alg_secantMethod} for completeness.

\begin{algorithm}
\caption{Secant-update method}
\label{alg_secantMethod}
\algsetup{indent=2em,linenodelimiter=}
\begin{algorithmic}[1]
\STATE{\text{\bf{Input:}} $b \in \mathbb{R}^m$,\, $y^{(0)}\in \mathbb{R}^{p}$, \, $\alpha^{(0)} >0$, \,$\eta > 1$, \,$\varepsilon \geq \|e\|$, and $\left\{U_p,\,B_p\right\}$ from the $p^{th}$ step of GKB}
\STATE{\text{\bf{Output:}} $\alpha^{(k+1)} \in \mathbb{R}$ }
\STATE{Compute $\tilde{r}_0^{(0)} = \tilde{r}_\alpha^{(0)}= U_{p}^Tb - B_py^{(0)} = \tilde{b}-B_py^{(0)}$}
\FOR{$k = 0,1,2,\dots$}
\STATE{$z=B_p^T\tilde{r}^{(k)}$}
\STATE{Compute $y_0^{(k+1)}=y_0^{(k)}+(B^T_{p}B_{p})^{-1}B_p^T\tilde{r}_0^{(k)}$}
\STATE{Compute $y_\alpha^{(k+1)} =y_\alpha^{(k)}+ \left(B^T_{p}B_{p}+\left[\alpha^{(k)} \right]^2I\right)^{-1}B_p^T\tilde{r}_\alpha^{(k)}$}
\STATE{$\tilde{r}_0^{(k+1)}= \tilde{b} - B_py_0^{(k+1)}$}
\STATE{$\tilde{r}_\alpha^{(k+1)}= \tilde{b} - B_py_\alpha^{(k+1)}$}
\STATE{$\gamma^{(k+1)}=\left\| \tilde{r}_0^{(k+1)} \right\|$}
\STATE{$\phi^{(k+1)}$}=$\left\| \tilde{r}_\alpha^{(k+1)} \right\|$
\STATE{$\alpha^{(k+1)} = \left| 
     \frac{\eta\varepsilon - \gamma^{(k+1)}}
          {\phi^{(k+1)} - \gamma^{(k+1)}}
   \right|\alpha^{(k)}$ }
\IF{$\phi^{(k+1)} \leq \eta \varepsilon$}
    \STATE{break}
\ENDIF
\ENDFOR
\end{algorithmic}
\end{algorithm}

\subsection{Nonstationary PIT} \label{sec:nonStationaryPIT}

We now turn our attention to the nonstationary PIT method whose iterates are computed possibly in a precision other than double. We first define our precision convention for computations and then unify Sections \ref{sec:stationaryPIT} and \ref{sec:regPar} to formalize the nonstationary PIT method. We close with a derivation and discussion of the nonstationary filter factors of PIT and comment on their computation in a precision possibly other than double.

When considering the computation of PIT in a precision possibly other than in fp64, we match the notational convention of the recent work \cite{NO24} and parallel the convention of recent works, e.g., \cite{Higham1,Higham22}, by using $\text{Pr}$ to denote the working precision used. In Section \ref{sec:results} we define the precisions we consider for our numerical experiments. We clarify that for a given precision used that all parts of the PIT algorithm are computed in this precision; this includes the GKB process, the secant-update method, and the PIT iterates themselves. That is, Pr is used uniformly in the PIT method.

After computing $p$ steps of the GKB process in precision Pr, the nonstationary PIT method may be represented at the $k^{th}$ iterate by
\begin{equation} \label{nonStationary_preconProj_Land}
    y^{(k)} = y^{(k-1)} + \left(B_p^TB_p + \left[\alpha^{(k)}\right]^2 I\right)^{-1}B^T_p\left(\tilde{b}-B_py^{(k-1)}\right).
\end{equation}
whose corresponding $k^{th}$ solution in $\mathbb{R}^n$ may be computed by $x^{(k)} = V_py^{(k)}$. In our numerical experiments we assume that the initial solution, i.e., $y^{(0)}$, is the zero vector so that the first projected iterate may be given by
\begin{equation*}
    \begin{split}
        y^{(1)} &= y^{(0)} + \left(B_p^TB_p + \left[\alpha^{(1)}\right]^2 I\right)^{-1}B^T_p\left(\tilde{b}-B_py^{(0)}\right)\\
        &= \left(B_p^TB_p + \left[\alpha^{(1)}\right]^2 I\right)^{-1}B^T_p\tilde{b}.
    \end{split}
\end{equation*}
Additionally, we assume $\alpha^{(0)} = 1$. This initial choice affects the number of iterations required by the PIT method to terminate, with a larger initial $\alpha$ typically requiring more iterations. We emphasize, however, that we have not experienced any detrimental effects due to variety of this choice. As is noted in Section \ref{sec:regPar}, the PIT method is terminated according to the discrepancy principle for the first iterate $k$ such that $\left\|\tilde{b}-B_py^{(k)}\right\|\leq \eta\varepsilon$. The nonstationary PIT method is provided in Algorithm \ref{alg1}. Henceforth, we drop explicit reference to the method being nonstationary and simply refer to the algorithm as PIT.

\begin{algorithm}
\caption{Nonstationary PIT method}
\label{alg1}
\algsetup{indent=2em,linenodelimiter=}
\begin{algorithmic}[1]
\STATE{\text{\bf{Input:}} $A \in \mathbb{R}^{m \times n}$, $b \in \mathbb{R}^m$, $x^{(0)} =0 \in \mathbb{R}^n$, $\alpha^{(1)}=1$, $\eta>1$, $\varepsilon \in \mathbb{R}$, and Pr}
\STATE{\text{\bf{Output:}} $x^{(k)} \in \mathbb{R}^n$}
\STATE{Compute $p$ steps of GKB: $AV_p = U_{p}B_{p}$}
\STATE{Compute $\tilde{r}^{(0)} = U_p^Tb = \tilde{b}$}
\FOR{$k = 1,2,\dots$}
\STATE{Compute $h^{(k-1)} = \left(B^T_{p}B_{p}+\left[\alpha^{(k)}\right]^2I\right)^{-1}B^T_{p}\tilde{r}^{(k-1)}$}
\STATE{Compute $y^{(k)} = y^{(k-1)} + h^{(k-1)}$}
\IF{$\left\|\tilde{b}-B_py^{(k)}\right\|\leq \eta\varepsilon$}
    \STATE{$x^{(k)} = V_py^{(k)}$}
    \STATE{break}
\ENDIF
\STATE{Compute $\tilde{r}^{(k)} = \tilde{b} - B_py^{(k)}$}
\STATE{Compute $\alpha^{(k)}$ using Algorithm \ref{alg_secantMethod}}
\ENDFOR
\end{algorithmic}
\end{algorithm}

\subsection{Filtered Solution of PIT} \label{sec:nonsta_filteredSoln}

To understand the filtering behavior of PIT when the regularization parameter is not fixed, we derive the filtered form of PIT recursively in Proposition \ref{prop:nonStationary_FF} so as to readily identify its filter factors. We note in Corollary \ref{coro:filtFactors} that when $\alpha^{(k)}$ are fixed for all $k$, that the result is equivalent to that of Proposition \ref{prop1}. We state this without proof since the result can also be shown in a straightforward manner by induction; for a proof when $M$ is a Tikhonov-type preconditioner, but not necessarily a low rank approximation of $A$ using a Krylov method, see Proposition 3.1 of \cite{NO24}.

\begin{prop}{(Nonstationary Filtered Solution of PIT)\\} \label{prop:nonStationary_FF}
After $p$ steps of GKB, the $k^{th}$ projected iterate, $y^{(k)}$, of PIT given by \eqref{nonStationary_preconProj_Land} may be written as 
\begin{equation} \label{nonstationary_PITA_projSoln}
	y^{(k)} = \sum_{i=1}^p \psi^{(k)}_i\frac{\hat{u}_{i}^T\tilde{b}}{\hat{\sigma}_i}\hat{v}_{i}, \quad \psi_i^{(k)}=\frac{\hat{\sigma}_i^2}{\hat{\sigma}_i^2 + \left[\alpha^{(k)}\right]^2}+\left(1 - \frac{\hat{\sigma}_i^2}{\hat{\sigma}_i^2 + \left[\alpha^{(k)}\right]^2} \right)\psi_i^{(k-1)}
\end{equation}
for $i=1,2,\dots,p$.
\begin{proof}
We first note that the $k^{th}$ iterate of PIT may be represented by
\begin{equation} \label{preconProj_Land_nonStationary}
\begin{split}
    y^{(k)} &= y^{(k-1)} + \left(B_p^TB_p + \left[\alpha^{(k)} \right]^2I\right)^{-1}B^T_p\left(\tilde{b}-B_py^{(k-1)}\right)\\
    &= y^{(k-1)} + \left(M_k^TM_k\right)^{-1}B^T_p\tilde{r}^{(k-1)}
\end{split}
\end{equation}
where we have used the definition of $M_k = \left(\hat{\Sigma}^T\hat{\Sigma} + \left[\alpha^{(k)}\right]^2I\right)^{1/2}\hat{V}^T = \hat{D}^{(k)}\hat{V}^T \in \mathbb{R}^{p \times p}$ and the SVD of $B_{p}=\hat{U}\hat{\Sigma}\hat{V}^T$ as previously used. We proceed via induction defining that $y^{(0)} =  0 \in \mathbb{R}^p$ so that there is nothing to show for the determination of $\psi_i^{(0)}=0$ for $i=1,\dots,p$.\\ \\
\emph{Base case $(k=1)$:} \\
Starting from the base relation given by \eqref{preconProj_Land_nonStationary} we have that
\begin{equation*}
    \begin{split}
        y^{(1)} &= y^{(0)} + \left(M_1^TM_1\right)^{-1}B_p^T\tilde{r}^{(0)} 
        = \left(M_1^TM_1\right)^{-1}B_p^T\tilde{b} \\
        &= \left(\hat{V}\left(\hat{D}^{(1)}\right)^2\hat{V}^T\right)^{-1}\hat{V}\hat{\Sigma}^T\hat{U}^T\tilde{b}.
    \end{split}
\end{equation*}
For notational simplicity we generically define $\bar{D}^{(k)} = \left(\hat{D}^{(k)}\right)^2$. Continuing, we may write
\begin{equation*}
    \begin{split}
        y^{(1)} &= \left(\hat{V}\bar{D}^{(1)}\hat{V}^T\right)^{-1}\hat{V}\hat{\Sigma}^T\hat{U}^T\tilde{b} \\
        &= \hat{V}\left(\bar{D}^{(1)}\right)^{-1}\hat{\Sigma}^T\hat{U}^T\tilde{b}\\ 
        &= \hat{V}\left(\bar{D}^{(1)}\right)^{-1}\hat{\Sigma}^T\hat{\Sigma}\hat{\Sigma}^{-1}\hat{U}^T\tilde{b}.
    \end{split}
\end{equation*}
Here we define the diagonal matrix $\Psi^{(1)} = \left(\bar{D}^{(1)}\right)^{-1}\hat{\Sigma}^T\hat{\Sigma} = \left(\bar{D}^{(1)}\right)^{-1}\hat{\Sigma}^2$ whose $i^{th}$ diagonal entry is given by $\psi^{(1)}_i = \hat{\sigma}_i^2/\left(\hat{\sigma}_i^2 + \left[\alpha^{(1)}\right]^2\right)$. With this we may continue to write
\begin{equation*}
    \begin{split}
        y^{(1)}
        &= \hat{V}\Psi^{(1)}\hat{\Sigma}^{-1}\hat{U}^T\tilde{b}\\
        &= \sum_{i=1}^p \left(\frac{\hat{\sigma}^2_i}{\hat{\sigma}_i^2+\left[\alpha^{(1)}\right]^2}\right)\frac{\hat{u}^T\tilde{b}}{\hat{\sigma}_i}\hat{v}_i\\
        &= \sum_{i=1}^p \left(\frac{\hat{\sigma}^2_i}{\hat{\sigma}_i^2+\left[\alpha^{(1)}\right]^2} + \left(\frac{\left[\alpha^{(1)}\right]^2}{\hat{\sigma}_i^2 + \left[\alpha^{(1)}\right]^2} \right)\psi_i^{(0)}\right)\frac{\hat{u}^T\tilde{b}}{\hat{\sigma}_i}\hat{v}_i\\
    \end{split}
\end{equation*}
completing the base case.\\ \\
\emph{Inductive step:} \\
Assume the inductive hypothesis is true for $k=n$ for $n>1$; that is
\begin{equation*}
    \begin{split}
        y^{(n)} &= \hat{V}\Psi^{(n)}\hat{\Sigma}^{-1}\hat{U}^T\tilde{b}\\
        &= \sum_{i=1}^p \left(\frac{\hat{\sigma}^2_i}{\hat{\sigma}_i^2+\left[\alpha^{(n)}\right]^2} + \left(\frac{\left[\alpha^{(n)}\right]^2}{\hat{\sigma}_i^2 + \left[\alpha^{(n)}\right]^2} \right)\psi_i^{(n-1)}\right)\frac{\hat{u}^T\tilde{b}}{\hat{\sigma}_i}\hat{v}_i
    \end{split}
\end{equation*}
so that the $i^{th}$ entry of the diagonal filtering matrix $\Psi^{(n)}$ is given by $\psi^{(n)}_i$. We proceed as follows
\begin{equation*}
    \begin{split}
        y^{(n+1)} &= y^{(n)} + \left(M_{n+1}^TM_{n+1}\right)^{-1}B_p^T\left(\tilde{b}-B_py^{(n)}\right)\\
        &= y^{(n)} + \left(M_{n+1}^TM_{n+1}\right)^{-1}B_p^T\tilde{b} + \left(M_{n+1}^TM_{n+1}\right)^{-1}B_p^TB_py^{(n)}\\
        &= \left(M_{n+1}^TM_{n+1}\right)^{-1}B_p^T\tilde{b} + \left(I - \left(M_{n+1}^TM_{n+1}\right)^{-1}B_p^TB_p\right)y^{(n)}\\
        &= \left(\hat{V}\bar{D}^{(n+1)}\hat{V}^T\right)^{-1}\hat{V}\hat{\Sigma}^T\hat{U}^T\tilde{b}+\left(I - \left(\hat{V}\bar{D}^{(n+1)}\hat{V}^T\right)^{-1}\hat{V}\hat{\Sigma}^T\hat{U}^T\hat{U}\hat{\Sigma}\hat{V}^T\right)y^{(n)}\\
        &= \hat{V}\left(\bar{D}^{(n+1)}\right)^{-1}\hat{\Sigma}^2\Sigma^{-1}\hat{U}^T\tilde{b} + \hat{V}\left(I-\left(\bar{D}^{(n+1)}\right)^{-1}\hat{\Sigma}^2\right)\hat{V}^Ty^{(n)}\\
        &= \hat{V}\Psi^{(n+1)}\Sigma^{-1}\hat{U}^T\tilde{b} + \hat{V}\left(I-\Psi^{(n+1)}\right)\hat{V}^T\hat{V}\Psi^{(n)}\hat{\Sigma}^{-1}\hat{U}^T\tilde{b}\\
        &= \hat{V}\left[\Psi^{(n+1)} + \left(I - \Psi^{(n+1)}\right)\Psi^{(n)}\right]\hat{\Sigma}^{-1}\hat{U}^T\tilde{b}\\
        &= \sum_{i=1}^p\left(\frac{\hat{\sigma}_i^2}{\hat{\sigma}_i^2 + \left[\alpha^{(n+1)}\right]^2}+\left(1 - \frac{\hat{\sigma}_i^2}{\hat{\sigma}_i^2 + \left[\alpha^{(n+1)}\right]^2} \right)\psi_i^{(n)}\right)\frac{\hat{u}^T\tilde{b}}{\hat{\sigma}_i}\hat{v}_i.
    \end{split}
\end{equation*}
Since both the base case and the inductive step have been shown, by mathematical induction the relation \eqref{nonstationary_PITA_projSoln} holds for every natural number $k$.
\end{proof}
\end{prop}

\begin{coro}{(Equivalence of Proposition \ref{prop1} and Proposition \ref{prop:nonStationary_FF} when $\alpha$ fixed)\\} \label{coro:filtFactors}
   If $\alpha^{(k)}>0$ is fixed as $\alpha$ for all $k=1,2,...$ then the recursively defined filter factors of Proposition \ref{prop:nonStationary_FF} coincide with those of Proposition \ref{prop1} when considering the same GKB subspace of size $p$, i.e.,
   $$y^{(k)} = \sum_{i=1}^p \psi^{(k)}_i\frac{\hat{u}_{i}^T\tilde{b}}{\hat{\sigma}_i}\hat{v}_{i}, \qquad \psi^{(k)}_i=1 - \left(\frac{\alpha^2}{\hat{\sigma_i}^2+\alpha^2}\right)^{k}$$
for $i=1,2,\dots,p$.
\end{coro}

We comment that if the projected residual $\left\|\tilde{b}-B_py^{(k)}\right\|$ is small after $k$ iterations, then the discrepancy principle \eqref{eq:secantDP} may be approximately satisfied. In this situation the next regularization parameter produced by Algorithm \ref{alg_secantMethod} will at least be as large as the previous one so that entry-wise change from the filter factors $\Psi^{(k-1)}$ to $\Psi^{(k)}$ is small. This can be understood by considering the $i^{th}$ filter factor represented in \eqref{nonstationary_PITA_projSoln}:
$$\psi_i^{(k)}=\frac{\hat{\sigma}_i^2}{\hat{\sigma}_i^2 + \left[\alpha^{(k)}\right]^2}+\left(1 - \frac{\hat{\sigma}_i^2}{\hat{\sigma}_i^2 + \left[\alpha^{(k)}\right]^2} \right)\psi_i^{(k-1)}.$$
Here, we may observe that if, for example, $\alpha^{(k)}>>\hat{\sigma}_i$, then the change between $\psi_i^{(k)}$ and $\psi_i^{(k-1)}$ for $i=1,2,\dots,p$ is expected to be small.

\section{Numerical Results and Preliminaries} \label{sec:results}

In this section we provide numerical experiments to illustrate the effectiveness of PIT computed in low precision compared to fp64 on the least-squares problem \eqref{minProb} and provide an exposition of the nonstationary filter factor result from Section \ref{sec:nonsta_filteredSoln}. The section is organized as follows: Sections \ref{secNum2} and \ref{secNum3} provide the background for the determination of the experimentally derived filter factors of PIT which we will term \emph{effective filter factors} of PIT on the least-squares problem and an overview of the simulation of low precision computations in MATLAB, respectively. In Section \ref{sec:experiments} we discuss our numerical results and their associated preliminaries.

\subsection{Effective filter factors} \label{secNum2}

To experimentally confirm our nonstationary filter factor result from Proposition \ref{prop:nonStationary_FF}, we compare against the corresponding effective filter factors. Using \eqref{filteredSoln}, the effective filter factors $\omega^{(k)}_j$ of the $k^{th}$ approximate projected solution of PIT using Algorithm \ref{alg1} may be computed as
\begin{equation} \label{effectiveFF}
    \omega^{(k)}_j = \frac{\hat{v}_{j}^Ty^{(k)}}{\hat{u}^T_{j}\tilde{b}}\hat{\sigma}_{j}, \quad j=1,2,\dots,n
\end{equation}
where we denote the vector containing the filter factors of the $k^{th}$ iterate by $\Omega^{(k)} \in \mathbb{R}^p$. We comment that when the $k^{th}$ iterate of PIT is compute using Algorithm \ref{alg1} using precision Pr that we also compute $\Omega^{(k)}$ in precision Pr using the numerical precision simulation software \textbf{chop}.

\subsection{Low precision simulation in MATLAB} \label{secNum3}

To consider the effectiveness of PIT in low precision and the experimental validity of our derived filter factors presented in Proposition \ref{prop:nonStationary_FF} we utilize the software package \textbf{chop} introduced by Higham and Pranesh in \cite{HighamPranesh19}. The \textbf{chop} function simulates lower precision arithmetic by rounding array entries given in a MATLAB native precision (e.g., single or double) to a target precision which is stored in a higher precision with non-utilized representation bits set to zero. The application of \textbf{chop} to structured discrete inverse problems was considered in \cite{Nagy23} with a focus on efficiency of use.

The \textbf{chop} software supports various common target precisions (e.g., single, half) as well as user-customizable ones. The precision choices and summary information we utilize in our numerical experiments in Section \ref{sec:experiments} are provided in Table \ref{choptable}. We comment that when accounting for orthogonality sensitivity in our GKB iterations that we did not experience any significant experimental differences when using precisions that did not support subnormal numbers, e.g., Google's bfloat16, compared to those based on the IEEE standard during our numerical investigations and therefore do not include them in our results.

\begin{table}[ht]
	\centering
    \caption{Precision choices considered in the numerical results: precision name, bits given for the mantissa and exponent, and the approximate decimal round-off error.}
    \begin{tabular}{c | c | c | c }
        \makecell{Pr\\Name} & \makecell{Exponent\\Bits} & \makecell{Mantissa\\Bits} & \makecell{Approx. Decimal\\Round-off} \\
        \ChangeRT{1.5pt}
        fp64 & $11$ & $52$ & $1.11e\text{-}16$ \\
        fp32 & $8$ & $23$ & $5.96e\text{-}8$ \\
        fp16 & $5$ & $10$ & $4.88e\text{-}4$ \\
    \end{tabular}
    \label{choptable}
\end{table}

The reported numerical results of this work were carried out in MATLAB R2024b 64-bit on a MacBook Pro laptop running MacOS Ventura with an Apple M2 Pro processor with @3.49 GHz and 16 GB of RAM. The computations other than those utilizing \textbf{chop} were carried out with about 15 significant decimal digits. SVD computations used for computing the theoretical filter factors from Proposition \ref{prop:nonStationary_FF} in fp16 were performed in Julia using the \emph{GenericLinearAlgebra} package.

\subsection{Examples and Results} \label{sec:experiments}

To effectively consider the filtering behavior of PIT and the behavior of the secant method for determining the regularization parameters at the $k^{th}$ iterate, we consider the 1D signal restoration problem \emph{Spectra} whose matrix models a symmetric Gaussian blur and $x$ is a simulated X-ray spectrum \cite{trussell}. The $a_{i,j}$ entries of $A\in \mathbb{R}^{64 \times 64}$ are given by
\begin{equation*} \label{spectraEntries}
    a_{i,j} = \frac{1}{\rho\sqrt{2\pi}}\exp\left(-\frac{(i-j)^2}{2\rho^2}\right),
\end{equation*}
with $\rho=2$ which results in a Toeplitz matrix. To realistically simulate the inverse problem, noise was added to the true right-hand side, $\hat{b}$, by forming the vector $e$ with normally distributed random entries with mean zero so that $b=Ax + e$; the vector $e$ is scaled so as to correspond to a specific noise that we will refer to as $\mu$, which is given by $\mu = 100\left(\|e\|/\|\hat{b}\|\right)$. The condition number of the matrix $A$ as determined by the MATLAB function $\text{cond}(\,)$ is $\approx10^9$; it can also be easily verified that the singular values of this matrix decay without a significant gap.

\begin{figure}[ht]
\begin{minipage}{1\linewidth}
	\centering
	\begin{minipage}{0.58\linewidth}
		\centering
		\includegraphics[width=\linewidth]{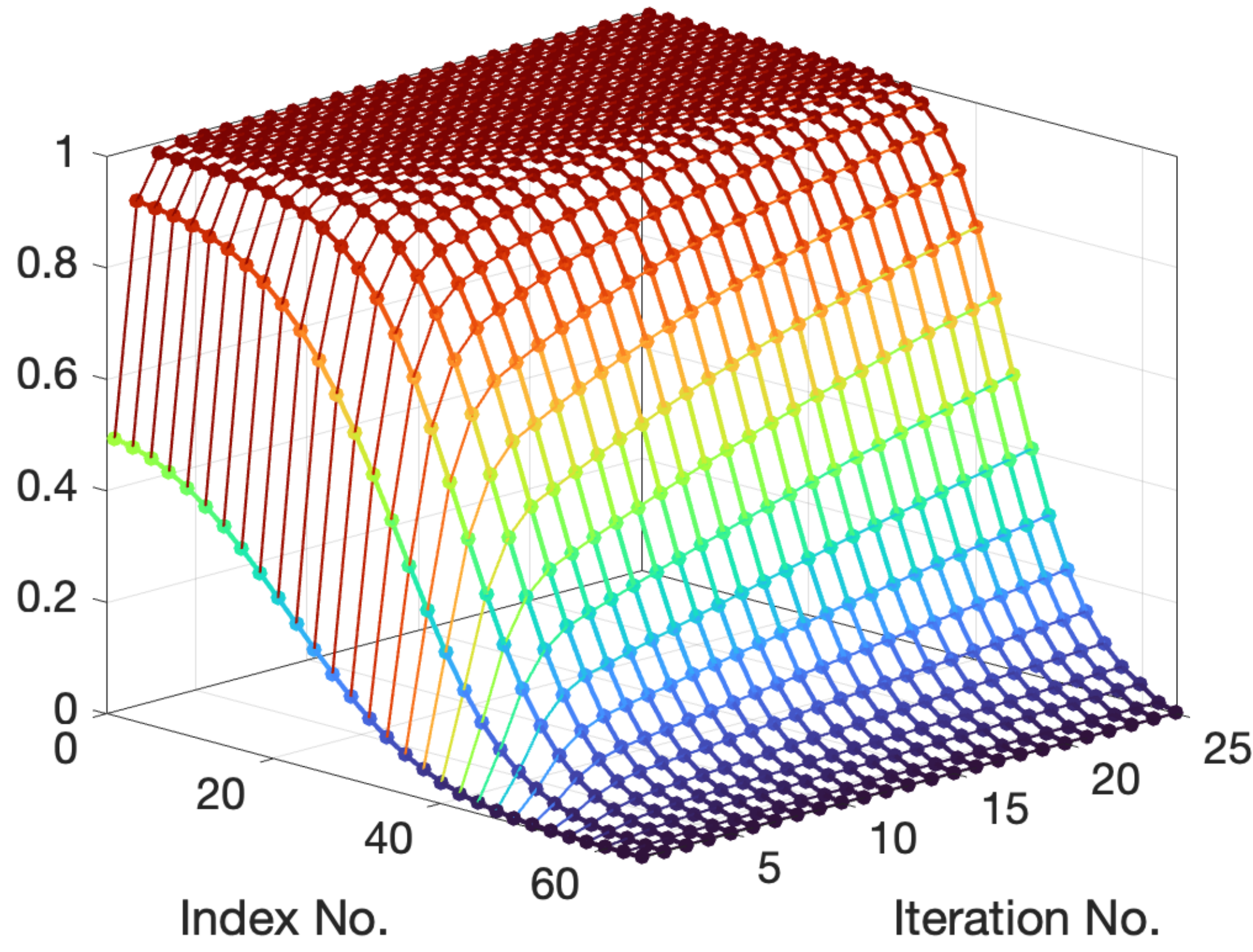}
	\end{minipage}
    \begin{minipage}{0.1\linewidth}
		\centering
		\includegraphics[width=0.85\linewidth]{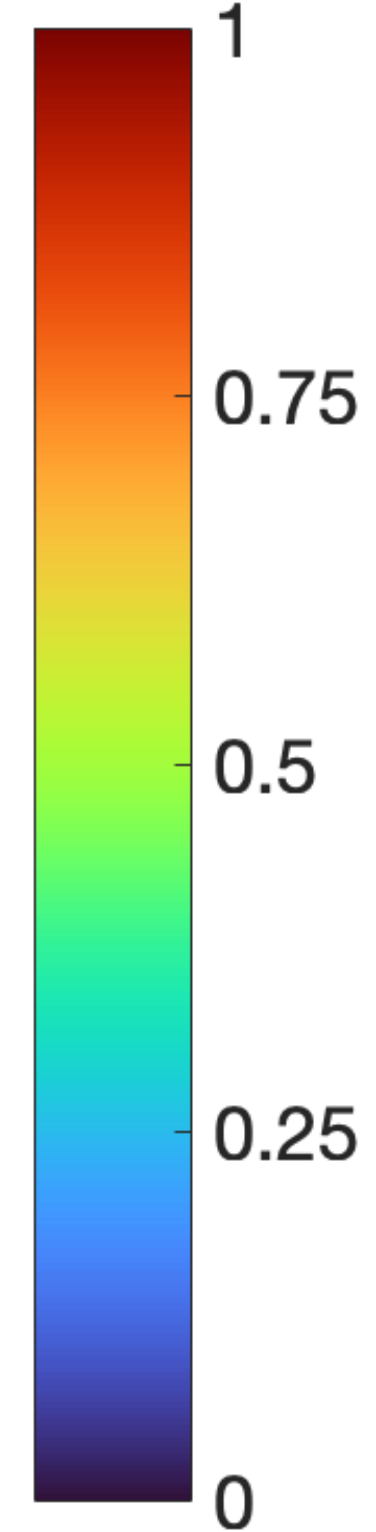}
	\end{minipage}
\end{minipage}
\caption{\emph{Spectra} example: effective filter factors for $25$ iterations of PIT in fp64 with a subspace of size $p=30$ and with $3\%$ noise. Each iteration represents $64$ filter factors whose values are denoted by $\omega_i^{(k)}$ and whose values may be read by the colorbar to the right or the value of the vertical axis.}
\label{fig:filterFactors}
\end{figure}

(\textbf{\emph{Spectra - filter factors and nonstationary regularization}}) - We first consider the filtering behavior of the PIT method with the secant method on the \emph{Spectra} problem contaminated by $3\%$ noise in fp64. For experimental focus, we compute a fixed number of steps, $p=30$, of GKB using full reorthogonalization for all three precisions considered. We consider the consequences that can occur from the loss of orthogonality of Krylov bases when computing in low precision in the image reconstruction examples to follow. For this experiment we computed the effective filter factors \eqref{effectiveFF} over 25 iterations and show these in Figure \ref{fig:filterFactors}. 

We may notice first and foremost that the filter factors of each experiment behave in a Tikhonov-type filtering manner - that is, early filter factors cluster near one and decay towards zero with increasing index. Since the regularization parameter is nonstationary for PIT, the secant method updates the regularization parameter $\alpha^{(k)}$ at each iteration $k$. We may observe visually in Figure \ref{fig:filterFactors} that the rate of change of the filter factors entry to entry (i.e., $\psi^{(k)}_i$ to $\psi^{(k+1)}_i$) is largest only for the first few iterations. This indicates that the iteration to iteration improvement in the solution becomes smaller as the iterations proceed. This occurs when the regularization parameter at iteration $k$ is computed to be larger than previous iterations causing the update to be minimal, which is computed by the secant method to coincide with the discrepancy principle being satisfied.

(\textbf{\emph{Spectra - filter factor comparison}}) - To confirm our derived filter factors of PIT from Proposition \ref{prop:nonStationary_FF} we compare against the effective filter factors that come from the same iteration $k$ of PIT computed by \eqref{effectiveFF}. We again utilize a noise level of $3\%$ and $30$ steps of GKB with full reorthogonalization. In Table \ref{table_FFs} we provide the summary statistics of the absolute difference between $\Phi^{(k)}$ and $\Omega^{(k)}$ of the 64 entries computed. We may observe that the mean absolute differences between the effective filter factors and those computed by Proposition \ref{prop:nonStationary_FF} approximately match the unit round off of the associated precision for the iterations listed (see Table \ref{choptable}). We only show iterations $1$, $5$, and $25$ for conciseness. For completeness we also show the min, max, and standard deviations of the absolute differences. We close by commenting that a finite precision error analysis is expected to support these results and to help explain the impact of the nonstationary preconditioner on the error; such an analysis is beyond the scope of this work.

\begin{table}[ht]
        \footnotesize
        \caption{\emph{Spectra} example: summary statistics of the absolute difference between the effective filter factors defined in \eqref{effectiveFF} and those computed by Proposition \ref{prop:nonStationary_FF} given by $|\Phi^{(k)} - \Omega^{(k)}|$ for precisions fp64, fp32, and fp16 with a noise level of $3\%$ and subspace of size $p=30$. We clarify that zero entries denoted by $0^{\ast}$ denote an absolute difference value that underflowed in the corresponding precision to $0$.}
	\centering
    \begin{tabular}{ c | c | c | c | c | c }
             \multirow{2}{*}{\scriptsize{Pr Choice}} & \multirow{2}{*}{\footnotesize{Iter.}} & \multicolumn{4}{c}{\footnotesize{Summary Stats. of $|\Phi^{(k)} - \Omega^{(k)}|$}} \\
            \cline{3-6}
            & & \footnotesize{Mean Error} & \footnotesize{Min Error} & \footnotesize{Max Error} & \footnotesize{Std. Dev.} \\
        \ChangeRT{1.5pt}
        \multirow{3}{*}{fp64} & $1$& $1.4e\text{-}16$ & $0^{\ast}$ & $6.7e\text{-}16$ & $1.8e\text{-}16$ \\
        \cline{3-6}
          & $10$& $4.4e\text{-}15$ & $0^{\ast}$ & $2.3e\text{-}15$ & $5.1e\text{-}16$ \\
        \cline{3-6}
          & $25$& $5.1e\text{-}16$ & $0^{\ast}$ & $2.1e\text{-}15$ & $5.0e\text{-}16$ \\
         \hline 
          \hline
          \hline
          \multirow{3}{*}{fp32} & $1$& $9.1e\text{-}8$ & $0^{\ast}$ & $5.7e\text{-}7$ & $1.4e\text{-}7$ \\
        \cline{3-6}
          & $10$& $4.9e\text{-}7$ & $2.0e\text{-}10$ & $3.4e\text{-}6$ & $6.7e\text{-}7$ \\
        \cline{3-6}
          & $25$& $5.7e\text{-}7$ & $5.0e\text{-}10$ & $3.6e\text{-}6$ & $7.3e\text{-}7$ \\
         \hline 
          \hline
          \hline
        \multirow{3}{*}{fp16} & $1$& $7.1e\text{-}4$ & $0^{\ast}$ & $1.3e\text{-}2$ & $2.5e\text{-}3$ \\
        \cline{3-6}
          & $10$& $4.5e\text{-}3$ & $0^{\ast}$ & $8.5e\text{-}2$ & $1.6e\text{-}2$ \\
        \cline{3-6}
          & $25$& $4.9e\text{-}3$ & $0^{\ast}$ & $9.0e\text{-}2$ & $1.6e\text{-}2$ \\
         \hline 
		\end{tabular}
\label{table_FFs}
\end{table}

\begin{figure}[ht]
\centering
\begin{subfigure}{0.32\textwidth}
    \centering
    \includegraphics[width=\linewidth,height=4cm,keepaspectratio]{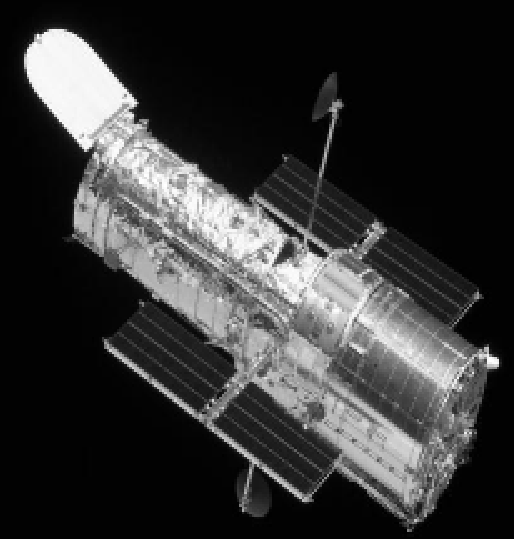}
    \caption{}
    \label{fig:xtrue}
\end{subfigure}
\hfill
\begin{subfigure}{0.32\textwidth}
    \centering
    \includegraphics[width=\linewidth,height=4cm,keepaspectratio]{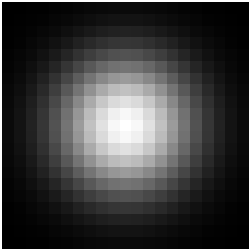}
    \caption{}
    \label{fig:prblur_psf}
\end{subfigure}
\hfill
\begin{subfigure}{0.32\textwidth}
    \centering
    \includegraphics[width=\linewidth,height=4cm,keepaspectratio]{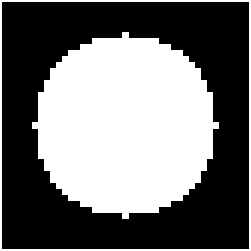}
    \caption{}
    \label{fig:prblurdef_psf}
\end{subfigure}
\vspace{0.75em}

\begin{subfigure}{0.32\textwidth}
\hfill
\end{subfigure}
\hfill
\begin{subfigure}{0.32\textwidth}
    \centering
    \includegraphics[width=\linewidth,height=4cm,keepaspectratio]{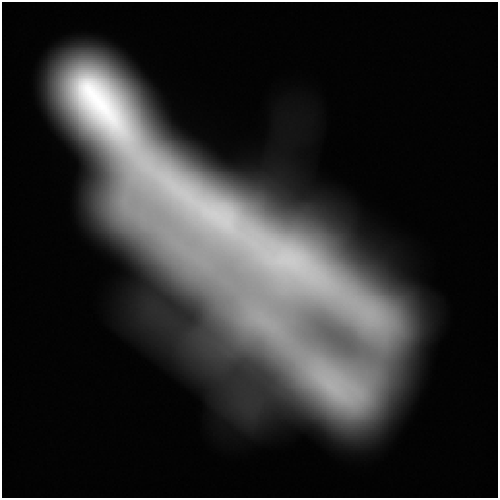}
    \caption{}
    \label{fig:prblur_bn}
\end{subfigure}
\hfill
\begin{subfigure}{0.32\textwidth}
    \centering
    \includegraphics[width=\linewidth,height=4cm,keepaspectratio]{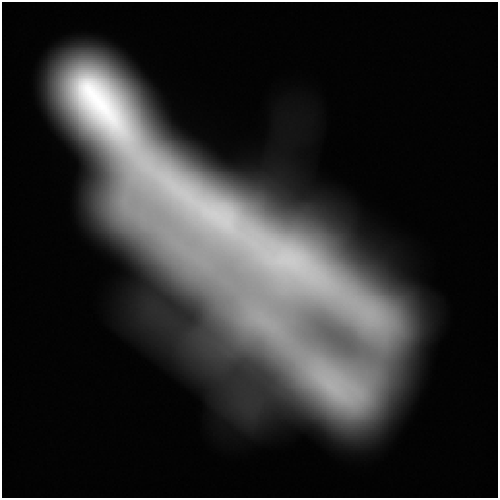}
    \caption{}
    \label{fig:prblurdef_bn}
\end{subfigure}
\hfill
\caption{2D image reconstruction examples: (A) true Hubble image ($256\times 256$ pixels), (B) \textit{Gauss} PSF, (C) \textit{Defocus} PSF, (D) Hubble image blurred by \textit{Gauss} with $1\%$ noise ($256\times 256$ pixels), and (E) Hubble image blurred by \textit{Defocus} with $1\%$ noise ($256\times 256$ pixels).}
\label{fig:test_problems}
\label{fig:layout}
\end{figure}

To evaluate the solution efficacy of PIT in low precision we consider two 2D image reconstruction problems from the IR Tools \cite{GHN18} software package: one a Gaussian blur and the other an out-of-focus blur. The point spread functions (PSFs) of both blurring operators are spatially invariant; see \cite{HNO} for more details. We will refer to these test problems as \textit{Gauss} and \textit{Defocus}, respectively. For both blurring operators, we use the Hubble space telescope image with a medium blur level and impose zero boundary conditions. The true image, as shown in Figure \ref{fig:xtrue}, contains $256 \times 256$ pixels. The corresponding blurring operators are of size $256^2 \times 256^2$ for each problem. The PSFs for \emph{Gauss} and \emph{Defocus} problems are shown in Figures \ref{fig:prblur_psf} and \ref{fig:prblurdef_psf}, respectively. 

The massive size of the operators make dense storage and application of \textbf{chop} to the full problems untenable. However, we may exploit the underlying Kronecker structure of these operators. For the \textit{Gauss} test problem, the PSF is rank-1 which allows for the blurring matrix $A$ to be written exactly as a single Kronecker product:
\[
A = A_r \otimes A_c 
  = (U_r\Sigma_rV_r^T)\otimes(U_c\Sigma_cV_c^T)
  = (U_r\otimes U_c)(\Sigma_r\otimes \Sigma_c)(V_r\otimes V_c)^T,
\]
where each factor is of a computationally manageable size ($256 \times 256$) to compute its SVD in the corresponding working precision. In contrast, the \textit{Defocus} PSF is not rank-1, so $A$ cannot be expressed as a single Kronecker product. Instead, it is approximated by a small sum of Kronecker products,
\[
A \approx A_K = \sum_{i=1}^{K} A_r^{(i)} \otimes A_c^{(i)},
\]
with $K=6$ terms providing an accurate representation for the $256\times256$ test problem \cite{nagy_mccormick_1998}. Each pair of Kronecker factors $A_r^{(i)}$ and $A_c^{(i)}$ is again small enough to allow for the efficient computation of their SVDs in the corresponding working precision. This allows us to use the same low-precision workflow as in the \emph{Gauss} case, at the cost of a factor of $K$ for additional storage and arithmetic. See \cite{Nagy23} for additional implementation details regarding the use of \textbf{chop} with Kronecker structures.

Similarly to the \textit{Spectra} example,  we contaminate the blurred images by Gaussian white noise with noise levels of $1\%$, $3\%$, and $5\%$. The noisy observed data for the $1\%$ case is shown in Figure \ref{fig:prblur_bn} and Figure \ref{fig:prblurdef_bn} for \textit{Gauss} and \textit{Defocus}, respectively. We implement PIT in precision Pr according to Algorithm \ref{alg1}, and employ the discrepancy principle with $\eta=1.01$ to determine the breakout iteration. To evaluate the quality of the reconstructed images, we compute the relative reconstructive error (RRE), which is defined by
$$\text{RRE}\left(x^{(k)}\right)=\frac{\left\| x^{(k)}-x\right\|}{\left\|x\right\|},$$
where $x^{(k)}$ denotes the computed approximate solution obtained when the breakout criterion is satisfied and $x$ represents the exact solution. In the numerical experiments that follow, we use the \textit{Gauss} test problem to investigate the effects of low-precision arithmetic on the orthogonality of the vectors generated by the GKB process during the PIT method. We then use the \textit{Defocus} test problem to focus on the choice of subspace size and efficacy of PIT in precisions lower than fp64.

\begin{figure}[ht]
    \centering
    \begin{subfigure}{0.32\textwidth}
        \centering
        \includegraphics[width=\linewidth]{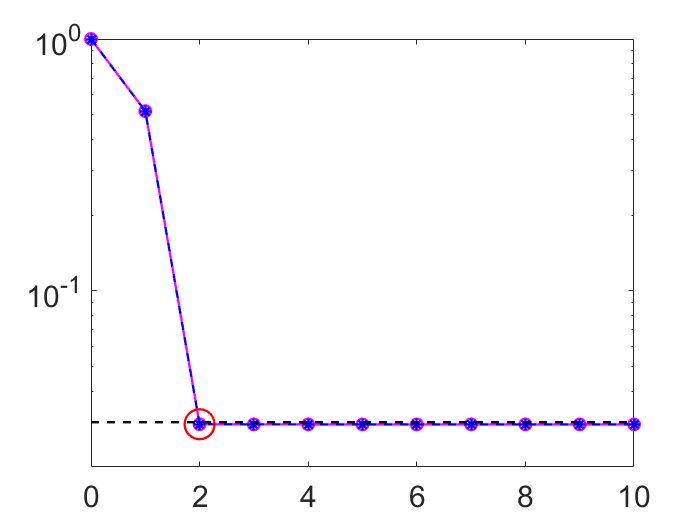}\\
        \vspace{0.5em}
        \includegraphics[width=\linewidth]{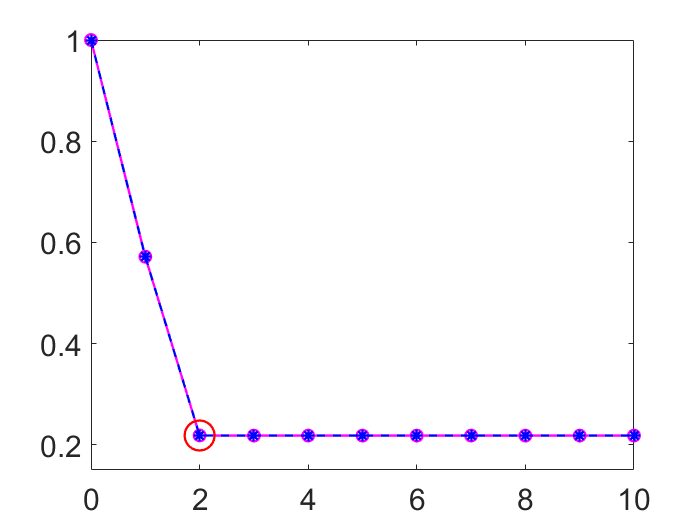}
        \caption{fp64}
    \end{subfigure}
    \hfill
    \begin{subfigure}{0.32\textwidth}
        \centering
        \includegraphics[width=\linewidth]{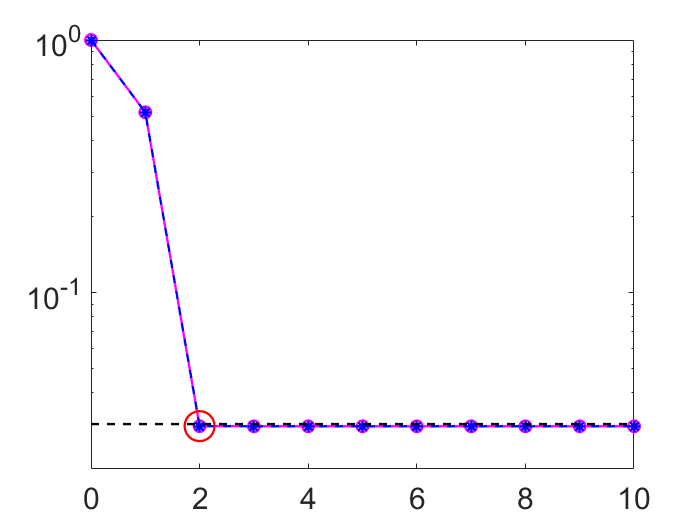}\\
        \vspace{0.5em}
        \includegraphics[width=\linewidth]{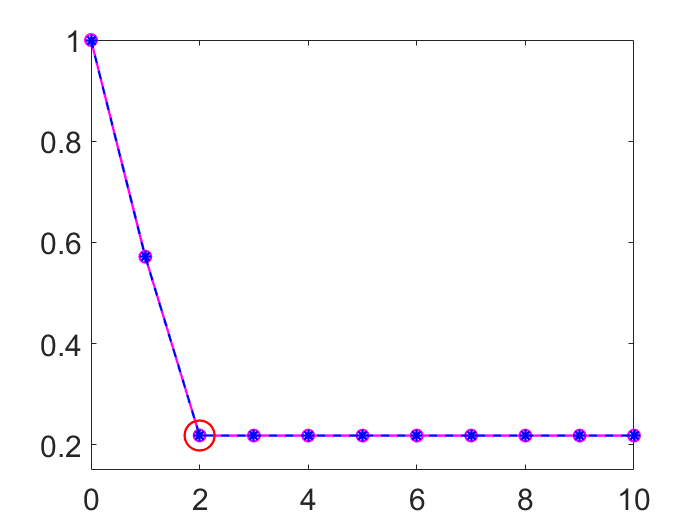}
        \caption{fp32}
    \end{subfigure}
    \hfill
    \begin{subfigure}{0.32\textwidth}
        \centering
        \includegraphics[width=\linewidth]{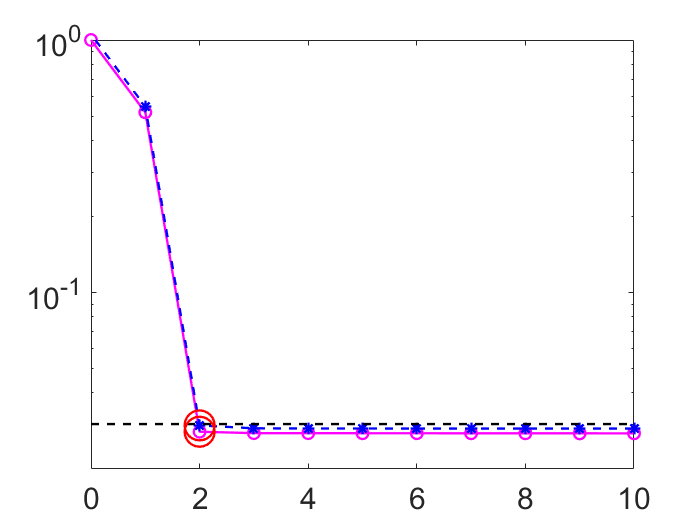}\\
        \vspace{0.5em}
        \includegraphics[width=\linewidth]{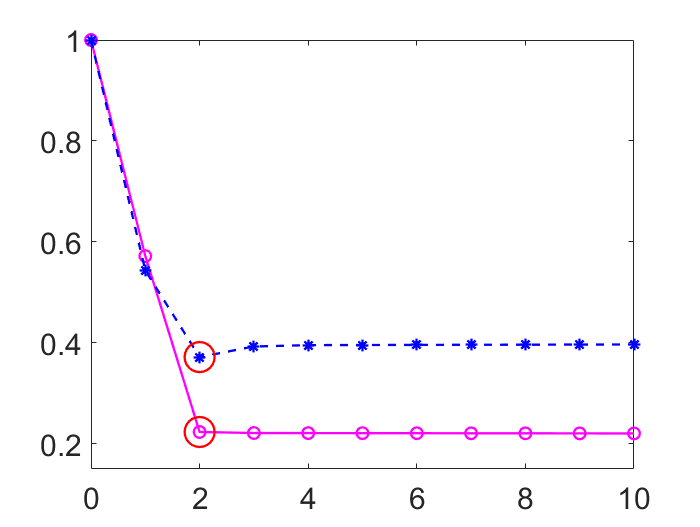}
        \caption{fp16}
    \end{subfigure}
    \caption{Relative residual (top row) and relative error (bottom row) against iteration number for the PIT method (with $30$ steps of GKB) applied to the \textit{Gauss} test case with $3\%$ noise, shown for the three precision levels: fp64, fp32, and fp16. The PIT method applied with reorthogonalization in the GKB process is represented by the solid magenta line and magenta circles. The PIT method without reorthogonalization is represented by the dashed blue line and blue asterisks. In the relative residual plots, the black dashed horizontal lines correspond to the discrepancy principle's breakout threshold. The red circles indicate the iteration at which the breakout criterion is satisfied.}
    \label{fig:Prblur_reorth}
\end{figure}

(\textbf{\emph{Gauss - orthogonality considerations for GKB}}) - We now consider the \textit{Gauss} test problem with $3\%$ noise. To evaluate the effects of enabling reorthogonalization (full reorthogonalization at each GKB step) during the PIT method, we present in Figure \ref{fig:Prblur_reorth} the relative residuals and RREs obtained by PIT with and without reorthogonalization for all three precisions considered. When reorthogonalization is disabled, we observe a significant loss in accuracy in fp16. However, across all three precision levels when reorthogonalization is enabled we observe that the RRE values for fp64, fp32, and fp16 are nearly indistinguishable, demonstrating that the dominant source of error in PIT without reorthogonalization is not the loss of precision within the projected iterations, but rather the stability of the Krylov basis used to form the projection. We observed similar results (not shown here) for $1\%$ and $5\%$ noise, and thus note that the use of GKB with reorthogonalization may be essential to achieve comparable accuracy when performing PIT in reduced precision.

(\textbf{\emph{Defocus - subspace size and solution efficacy}}) -
For our final numerical investigation we evaluate the behavior of PIT on the \textit{Defocus} test problem for both $1\%$ and $5\%$ noise levels for all three precisions considered. For this example, we enable reorthogonalization for all GKB steps. The numerical results are reported in Table \ref{tab:defocus}. Here, we consider both the number of GKB steps, $p$, and the iteration count when the discrepancy principle is satisfied, or when the maximum allowed iterations are computed; for the ladder we set this to be $20$. 

\begin{figure}[ht]
\centering
\begin{subfigure}{0.32\textwidth}
    \centering
    \includegraphics[width=\linewidth,height=4cm,keepaspectratio]{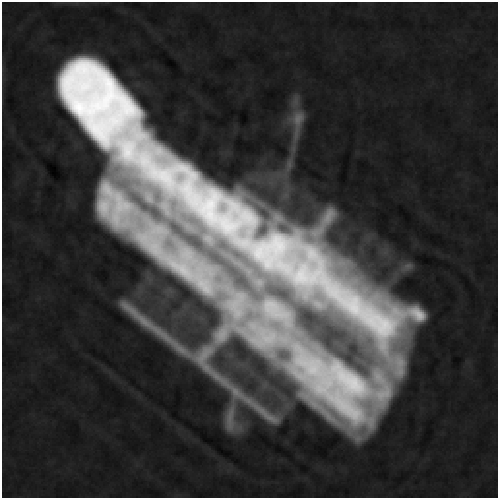}
    \caption{fp64}
\end{subfigure}
\hfill
\begin{subfigure}{0.32\textwidth}
    \centering
    \includegraphics[width=\linewidth,height=4cm,keepaspectratio]{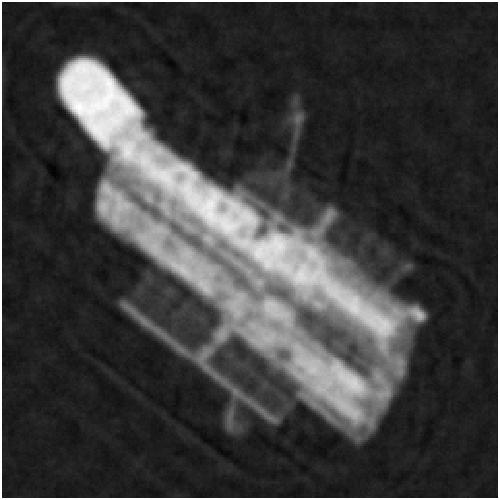}
    \caption{fp32}
\end{subfigure}
\hfill
\begin{subfigure}{0.32\textwidth}
    \centering
    \includegraphics[width=\linewidth,height=4cm,keepaspectratio]{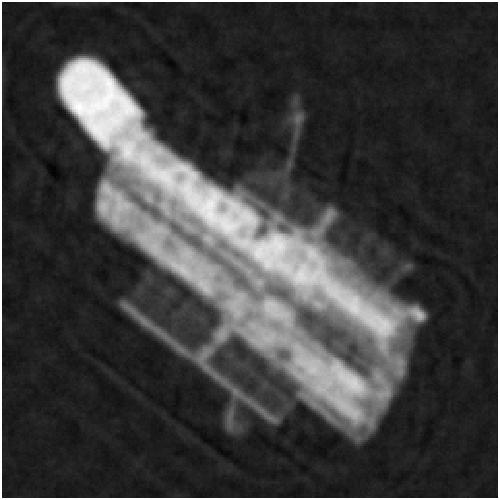}
    \caption{fp16}
\end{subfigure}
\hfill
\caption{2D image reconstructions for the PIT method applied to the \textit{Defocus} example with $1\%$ noise, $35$ GKB steps, and precision levels fp64, fp32, and fp16.}
\label{fig:defocus_reconstructions}
\end{figure}

We remark that the resulting RREs are consistent across all precision levels for every pairing of noise level and choice of subspace size $p$. This agrees with our earlier observations that, once a sufficiently stable (with respect to orthonormality of the basis vectors) Krylov subspace is formed, the projected iterations that PIT is computed in, for example fp16, do not appear to be negatively affected by the computations in low precision. Thus, similarly to the \emph{Gauss} example we observe that the solution efficacy appears to be dominated by the quality of the GKB basis rather than the arithmetic precision used within PIT. Furthermore, we comment that our studies here in precisions lower than fp64 align with the findings from \cite{BGOPR25} that a sufficiently large subspace is required to achieve breakout according to the discrepancy principle after a few iterations. Lastly, we point out in Figure \ref{fig:defocus_reconstructions} that the reconstructions using PIT in fp64, fp32, and fp16 look essentially indistinguishable to the naked eye which further supports the use of low precision for image deconvolution applications.

\begin{table}[ht]
\footnotesize
\caption{RRE results for the PIT method applied to the \textit{Defocus} example with $1\%$ and $5\%$ noise. The precision levels fp64, fp32, and fp16 are considered and the number of steps of GKB allotted during PIT is varied ($p=25,30,35$). The RRE values presented correspond to the $k^{th}$ iterative solution when the discrepancy principle is satisfied, or alternatively, when the specified maximum number of iterations ($k=20$) is achieved.}
\centering
\begin{tabular}{ c | c | c | c | c }
    \scriptsize{Noise} 
      & \scriptsize{\makecell{GKB Steps\\($p$)}}
      & \scriptsize{\makecell{Pr\\Choice}}
      & \scriptsize{RRE}
      & \scriptsize{Iter.} \\
\ChangeRT{1.5pt}
\multirow{9}{*}{1\%} 
    & \multirow{3}{*}{25} 
        & fp64 & 0.2186 & 20 \\   
        & & fp32 & 0.2186 & 20 \\
        & & fp16 & 0.2206 & 20 \\
    \cline{2-5}
    & \multirow{3}{*}{30} 
        & fp64 & 0.2097 & 20 \\
        & & fp32 & 0.2097 & 20 \\
        & & fp16 & 0.2102 & 20 \\
    \cline{2-5}
    & \multirow{3}{*}{35} 
        & fp64 & 0.2038 & 2 \\
        & & fp32 & 0.2038 & 2 \\
        & & fp16 & 0.2084 & 2 \\
\hline \hline \hline
\multirow{9}{*}{5\%} 
    & \multirow{3}{*}{25} 
        & fp64 & 0.2549 & 2 \\
        & & fp32 & 0.2549 & 2 \\
        & & fp16 & 0.2600 & 2 \\
    \cline{2-5}
    & \multirow{3}{*}{30} 
        & fp64 & 0.2572 & 3 \\
        & & fp32 & 0.2572 & 3 \\
        & & fp16 & 0.2641 & 2 \\
    \cline{2-5}
    & \multirow{3}{*}{35} 
        & fp64 & 0.2553 & 3 \\
        & & fp32 & 0.2553 & 3 \\
        & & fp16 & 0.2548 & 3 \\
\ChangeRT{0.4pt}
\end{tabular}
\label{tab:defocus}
\end{table}

\section{Conclusions} \label{sec:conc}
In this work, we investigated the regularizing behavior of the PIT method in low precision. To do so, we considered it through the lens of a filtering method, writing it as a projected preconditioned Landweber method with a Tikhonov-type preconditioner using the GKB process. In doing so we confirmed that the projected filter factors behave in a Tikhonov-type manner in either low (i.e., fp16) or high (i.e., fp64) precision and discussed their interaction with nonstationary regularization parameters. We also considered a secant-update method in conjunction with the discrepancy principle for the determination of the regularization parameters at each step of PIT, which was observed to be robust for all precisions considered herein. Our numerical results suggest that the basis vectors of the formed Krylov subspace can lose orthonormality when computed in sufficiently low precision, which we showed can lead to poor relative errors for image deconvolution problems. However, we found in our experiments that performing reorthogonalization can alleviate this issue and that PIT computed in low precision can deliver comparable working accuracy (i.e., to within a few decimal places in relative error) for image deconvolution problems with a sufficiently large subspace compared to computations carried out in double precision.  

\bibliographystyle{siam}
\bibliography{references}

\end{document}